\documentclass[final]{amsart}
\usepackage{amssymb}
\usepackage{graphicx}
\def\fig#1{
\includegraphics[height=3.5cm]{#1.eps}
}
\newcommand{\QQ}[2]{Q_{#1}(#2;t,x,s)}

\renewcommand{\boxtimes}{\bigstar_c}

\def\<{\langle}
\def\>{\rangle}

\def\<{\langle}
\def\>{\rangle}

\newcommand{\rf}[1]{(\ref{#1})}

\newcommand{\calF}{{\mathcal F}}
\newcommand{\calG}{{\mathcal G}}

\newcommand{\calL}{{\mathcal L}}

\newcommand{\aaa}{{\mathbf a}}
\newcommand{\bbb}{{\mathbf b}}

\newcommand{\AAA}{{\mathbf A}}
\newcommand{\BBB}{{\mathbf B}}
\newcommand{\CCC}{{\mathbf C}}
\newcommand{\DDD}{{\mathbf D}}

\newcommand{\sR}{{\mathbb R}}

\newcommand{\sN}{{\mathbb N}}

\newcommand{\Var}{\mbox{Var}}
\newcommand{\eps}{\varepsilon}

\newcommand{\be}{\begin{equation}}
\newcommand{\ee}{\end{equation}}

      \newtheorem{theorem}{Theorem}[section]
       \newtheorem{proposition}[theorem]{Proposition}
       \newtheorem{corollary}[theorem]{Corollary}
       \newtheorem{lemma}[theorem]{Lemma}

\theoremstyle{remark}
       \newtheorem{remark}{Remark}[section]
\theoremstyle{definition}

\newtheorem{assumption}{Assumption}

\def\<{\langle}
\def\>{\rangle}

\def\<{\langle}
\def\>{\rangle}

\author{W{\l}odzimierz  Bryc
}
\thanks{\noindent Research partially supported by NSF
grant \#INT-0332062 and by the C.P. Taft Memorial Fund.}
\address{
Department of Mathematics,
University of Cincinnati,
PO Box 210025,
Cincinnati, OH 45221--0025, USA}
\email{Wlodzimierz.Bryc@UC.edu}

\author{Jacek Weso{\l}owski}
\address{ Faculty of Mathematics and Information Science\\
Warsaw University of Technology\\ pl. Politechniki 1\\ 00-661
Warszawa, Poland}
\email{wesolo@alpha.mini.pw.edu.pl}

\date{March 8, 2004; Revised \today\ {\tt File: \jobname.TEX}}
\subjclass[2000]{60J25}

\title
{Bi-Poisson process}
\begin{document}
\begin{abstract}
We study a two parameter family of processes with  linear regressions and linear
conditional variances. We give conditions for the unique solution of this problem, and point  out the
connection between the resulting Markov processes and the  generalized convolutions introduced by
Bo\.zejko and Speicher  \cite{Bozejko-Speicher91}.
\end{abstract}
\maketitle

\section{Introduction}
Throughout this paper $(X_t)_{t\geq 0}$ is a 
square
integrable stochastic process such that for all $t,s\geq 0$
\begin{equation}\label{EQ: cov}
E(X_t)=0,\: E(X_tX_s)=\min\{t,s\}.
\end{equation}
Consider the $\sigma$-fields
 $\calG_{s,u}=\sigma\{X_t: t\in[0,s]\cup[u,\infty)\}$,
 $\calF_s=\sigma\{X_t: t\in[0,s]\}$, $\calG_u=\sigma\{X_t: t\in[u,\infty)\}$.
We  assume that the process has linear  regressions,
\begin{assumption} \label{ASSUME: 1} For all $0\leq s<t<u$,
\begin{equation}
\label{EQ: LR} E(X_t|\calG_{ s, u})=\aaa X_s+\bbb X_u,
\end{equation}
where
\begin{equation}\label{EQ: a+b}
\aaa=\aaa(t|s,u)=\frac{u-t}{u-s},\:
\bbb=\bbb(t|s,u)=\frac{t-s}{u-s}.
\end{equation}
 are the deterministic functions of $0\leq s<t<u$.
\end{assumption}

We also assume that the process has  quadratic
 conditional variances;
\begin{equation}
\label{EQ: QV} E(X_t^2|\calG_{s,u})
=\AAA X_s^2+\BBB X_sX_u+\CCC X_u^2+\DDD+\alpha X_s+\beta X_u,
\end{equation}
where
$\AAA =\AAA(t|s,u), \BBB =\BBB(t|s,u), \CCC =\CCC(t|s,u),
 \DDD=\DDD(t|s,u), \alpha=\alpha(t|s,u),
\beta=\beta(t|s,u)$ are the deterministic functions of $0<s<t<u$.

Generically, conditions \rf{EQ: cov},  \rf{EQ: LR},
and  \rf{EQ: QV} imply that there are
five real parameters $q,\eta,\theta,\sigma,\tau$ such that

\begin{eqnarray}
\label{EQ: A} \AAA(t|s,u)&=&\frac{(u-t)(u(1+\sigma t)+\tau-qt)}{(u-s)(u(1+\sigma s)+\tau-qs)},\\
\label{EQ: B}\BBB(t|s,u)&=&\frac{(u-t)(t-s)(1+q)}{(u-s)(u(1+\sigma s)+\tau-qs)},\\
\label{EQ: C}\CCC(t|s,u)&=&\frac{(t-s)(t(1+\sigma s)+\tau-qs)}{(u-s)(u(1+\sigma s)+\tau-qs)} ,\\
\label{EQ: D}\DDD(t|s,u)&=&\frac{(u-t)(t-s)}{u(1+\sigma s)+\tau-qs},\\
\alpha(t|s,u)  &=& \frac{(u-t)(t-s)}{u(1+\sigma s)+\tau-qs} \times \frac{u\eta-\theta}{u-s}, \label{EQ: alpha}
\\
\beta(t|s,u) &=& \frac{(u-t)(t-s)}{u(1+\sigma
s)+\tau-qs}\times\frac{\theta-s\eta}{ u-s}. \label{EQ: beta}
\end{eqnarray}
This gives after a calculation
\begin{eqnarray}\label{EQ: q-Var}
&\Var(X_t|\calG_{s,u })= \frac{(u-t)(t-s)}{u(1+\sigma
s)+\tau-qs}\left( 1+ \sigma \frac{(uX_s-sX_u)^2}{(u-s)^2}+\eta
\frac{uX_s-sX_u}{u-s} \right. &\\ \nonumber &\left.
+\tau\frac{(X_u-X_s)^2}{(u-s)^2}+\theta\frac{X_u-X_s}{u-s}
+(1-q)\frac{(X_u-X_s)(sX_u-uX_s)}{(u-s)^2} \right),&
\end{eqnarray}
compare \cite[Proposition 2.5]{Bryc-Wesolowski-03}. (Recall that
the conditional variance of $X$  with respect to a $\sigma$-field
$\calF$ is defined as
$\Var(X|\calF)=E(X^2|\calF)-\left(E(X|\calF)\right)^2$.)

 In \cite{Bryc-Wesolowski-03} we prove that the solution
of equations \rf{EQ: cov}, \rf{EQ: LR},  \rf{EQ: q-Var}  exists
and is unique when $-1<q\leq 1$, and $\sigma=\eta=0$; it is then given
by the Markov process which we called $q$-Meixner process. (The case $q=1$ yields L\'evy processes, and
was studied earlier by several authors, see \cite{Wesolowski93}, and the references therein.)
Due to
the invariance of this problem under the symmetry that maps
$(X_t)$ to the process $(tX_{1/t})$, processes that satisfy
\rf{EQ: q-Var} with $-1<q\leq 1, \tau=\theta=0$ are also Markov,
and can be expressed in terms of the $q$-Meixner processes as
$tX_{1/t}$. The main feature of these examples are trivial
(constant) conditional variances in one direction of time, which
leads to technical simplifications.

The study of the remaining
 cases poses difficulties, as several steps from
 \cite{Bryc-Wesolowski-03} break down.
 In this paper we consider the next simplest case, which one may call
 the free bi-Poisson processes.
The $q$-Poisson processes,
in particular, the classical Poisson process and the free Poisson process, have linear conditional
 variances when conditioned with respect to the future,
and constant conditional variances when conditioned with respect to
the past. The bi-Poisson process has linear conditional variances under each uni-directional
 conditioning; it corresponds to the choice of
 $\sigma=\tau=0$ in \rf{EQ: q-Var}.
The adjective "free" refers to $q=0$.  The role of these simplifying conditions seems technical: linear conditional variances
imply that all moments are finite, see Lemma
 \ref{Lemma:moments}; additional condition that $q=0$ allows us to guess useful
 algebraic identities between the orthogonal polynomials
 in Proposition \ref{THM: poly algebra}.
These considerations lead to the following.
\begin{assumption}\label{ASSUME: 2} For all $0\leq s<t<u$,
\begin{eqnarray}\label{EQ: 0-Var}
&\Var(X_t|\calG_{s,u})=&\\
\nonumber
&\frac{(u-t)(t-s)}{u}\left(
1+
\eta  \frac{uX_s-sX_u}{u-s} +\theta\frac{X_u-X_s}{u-s}
+\frac{(X_u-X_s)(sX_u-uX_s)}{(u-s)^2}
\right).&
\end{eqnarray}
\end{assumption}
In Section \ref{Sect:free Poisson}  we construct the Markov process
 with covariances \rf{EQ: cov}, linear regressions \rf{EQ: LR},
and conditional variances \rf{EQ: 0-Var} for a large set of real
 parameters $\eta,\theta$.
In Section \ref{Sect:Char} we show that the solution is unique. In Section \ref{Sect: gen con}
we point out that when $\theta=1$ the one-dimensional distributions of the bi-Poisson process
are closed under a generalized free convolution.

\section{Existence}\label{Sect:free Poisson}
If $\eta=0$, formula \rf{EQ: 0-Var} coincides with
\cite[(28)]{Bryc-Wesolowski-03} with $\tau=q=0$, so the
corresponding Markov process exists and is determined uniquely, see
\cite[Theorem 3.5]{Bryc-Wesolowski-03}. Since the transformation
$X_t\mapsto tX_{1/t}$ switches the roles of $\eta,\theta$, the case
$\theta=0$ follows, too. We may therefore restrict our attention
to the case $\eta\theta\ne 0$. The construction of the processes
is based on the idea already exploited in
\cite{Bryc-Wesolowski-03}; namely,
 we construct the transition probabilities of the suitable Markov process,
  by defining  the corresponding orthogonal polynomials.
Under current assumptions, this task requires more work as we need to ensure
that the coefficient at the third term of the recurrence for the polynomials  is non-negative.
The construction relies on  new identities between the
 orthogonal polynomials, which are used to verify the martingale
 polynomial property \rf{EQ: proj}; the latter property fails for more general
 values of parameters in \rf{EQ: q-Var}.

 \subsection{One dimensional distributions}
We begin by carefully examining the "candidate" for the one
 dimensional distribution of $X_t$.  For $t>0$,
let $p_0(x;t)=1$, and
consider the following monic polynomials $\{p_n(x;t):\ n\geq 1\}$ in
 variable $x$.

\begin{eqnarray}
xp_0&=&p_1+0p_0  \label{EQ: p0},\\
 xp_1&=&p_2+(t\eta+\theta)p_1+tp_0 \label{EQ: p1},\\
xp_n&=&p_{n+1}+(t\eta+\theta)p_n+t(1+\eta\theta)p_{n-1}, \, n\geq 2.
 \label{EQ: pn rec}
\end{eqnarray}

From the general theory of orthogonal polynomials, if $1+\eta\theta
 \geq 0$ then there exists a unique
probability measure $\pi_t$ such that $p_n(x;t)$ are orthogonal with
 respect to $\pi_t$, see \cite{Chihara}.
We will need the following.
\begin{lemma}\label{THM: support}
\begin{equation}\label{EQ: 1+x eta}
\pi_t(\{x: 1+\eta x <0\})=0.
\end{equation}
\end{lemma}
\begin{figure}[h]
\fig{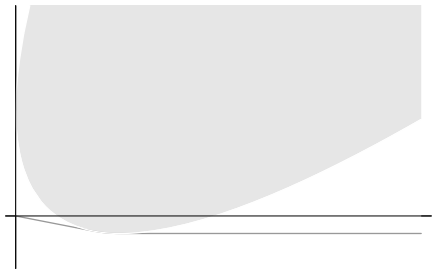}
\caption{The vertical cross-section of the gray area is the support
 of  $\pi_t$ for $t>0$;
the gray lines represent the support of
the  discrete part. This picture represents  the case $\eta>0,\theta>0$.\label{FIG: 1}}
\end{figure}


\begin{figure}[h]
\fig{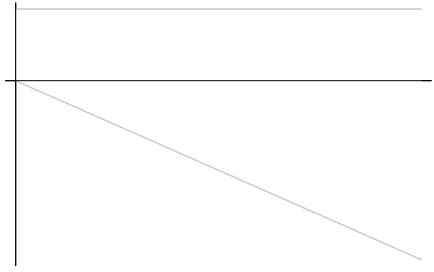}
\caption{The  support of  the distribution of $X_t$
for $t>0$ in the degenerate case $\eta<0,\theta>0,
 1+\eta\theta=0$.\label{FIG: 3}}
\end{figure}

\begin{proof}[Proof of Proposition \ref{THM: support}]
There is nothing to prove when $\eta=0$, so without loss of
 generality we assume that $\eta\ne0$.

If $1+\eta\theta=0$ then the recurrence is degenerate and the
 distribution is
supported at zeros of polynomial $p_2(x)=x^2-(t\eta+\theta)x-t$;
this follows from the fact that all higher order polynomials are
 multiples of $p_2$.
The support $\mbox{supp} (\pi_t)=\{-t/\theta,-1/\eta\}$, see Fig.
 \ref{FIG:  3}, is
disjoint with  the open set $\{x: 1+\eta x <0\}$,  ending the proof
 in this case.

If $1+\eta\theta>0$, then  \rf{EQ: pn rec} is a constant coefficient
 recurrence which has been analyzed by several
 authors, see \cite{Saitoh-Yoshida01}.
The Cauchy transform $$G(z)=\int \frac{1}{z-x}\pi_t(dx)$$ is given by
 the corresponding continued fraction,
$$
G(z)=\displaystyle\frac{1}
{z-\displaystyle\frac{t}
{z-(t\eta+\theta)-\displaystyle\frac{t(1+\eta\theta)}{z-(t\eta+\theta)-\displaystyle\frac{t(1+\eta\theta)}{\ddots}}}}
$$
which after a calculation gives
\begin{equation}\label{EQ: G(z)}
\displaystyle
G_{}(z)=\displaystyle
\frac{z(1+2\eta\theta)+t\eta+\theta-\sqrt{(z-(t\eta+\theta))^2-4 t(1+\eta\theta)}}
{2(1 + z \eta) (t + z \theta)}
.
\end{equation}
The Stieltjes inversion formula gives  the distribution $\pi_t$ as the limit in distribution
as $\eps\to 0^+$ of the absolutely continuous measures $-\frac{1}{\pi} \Im G(x+i\eps) dx$. This
 gives

\begin{equation}\label{EQ: distr}
\pi_{t}(dx)=\frac{t}{2\pi}\frac{
\sqrt{4 t (1+\eta\theta)-(x-t\eta-\theta)^2}}
{(x\eta+1)(x\theta+t)}
1_{(x-t\eta-\theta)^2<4t(1+\eta\theta)}
+ p(t) \delta_{-t/\theta}+q(t) \delta_{-1/\eta}.
\end{equation}
The weights at the discrete point masses are
$$
p(t)=\frac{-((1+\eta\theta)t-\theta^2)/\theta+\eps
|(1+\eta\theta)t-\theta^2|/|\theta|}{2(\theta-\eta t)}
$$
and
$$
q(t)=\frac{\eta(t-(1+\eta\theta)/\eta^2)+\eps|\eta||t-(1+\eta\theta)/\eta^2|}{2(\eta
t-\theta)},
$$
where the sign $\eps=\eps(t,\eta,\theta)=\pm 1$ is selected
 simultaneously for both expressions by the appropriate
 choice of the branch of the square root. We found that a
 practical way to choose the sign is to select
 $\eps=\pm1$ so that both  expressions give a number in the interval
 $[0,1]$; in our setting  this determines
 $\eps$ uniquely for every choice of parameters, after all the cases
 are considered.

It is easy to check that the support of the absolutely continuous part of $\pi_t$ does not intersects the set
$\{x:1+x\eta<
 0\}$. The support of the discrete part consists of at most two-points: $\{-t/\theta,\-1/\eta\}$.
  Thus the only possibility for the set $\{x:1+x\eta<
 0\}$ to carry positive $\pi_t$-probability is
 when $-t/\theta\in \{x:1+x\eta<0\}$. This is possible only if
 $\eta\theta>0$ and $t$ is large enough. The
Stieltjes inversion formula gives the weight of $-t/\theta$ as
 $$
 p(t)=\frac{\left(\theta^2-(1+\eta\theta)t\right)_+}{\theta^2-t\eta \theta }.$$
Thus the point $-t/\theta$ carries positive probability $p(t)$ only for
 $t<\frac{\theta^2}{1+\eta\theta}\leq \theta/\eta$;
on the other hand,  $-t/\theta\in \{x:1+x\eta<0\}$ only for $t>\theta/\eta$.

\end{proof}

\subsection{Transition probabilities}
 Fix $0<s<t$, and let $x\in\sR$ be such that $1+x\eta\geq0$.
We define monic polynomials in variable $y$  by the three-step
 recurrence
$$\QQ{0}{y}=1,$$
$$\QQ{1}{y}=y-x,$$
 $$
     y\QQ{1}{y}=\QQ{2}{y}+ ((t-s)\eta+\theta)\QQ{1}{y}+
     (t-s)(1+x\eta)\QQ{0}{y},
     $$
 and for $n\geq 2$ by the constant coefficients recurrence
   \begin{equation}\label{EQ: 3-step Q}
     y\QQ{n}{y}=\QQ{n+1}{y}+(t\eta+\theta)\QQ{n}{y}+t(1+\eta\theta)\QQ{n-1}{y}.
\end{equation}

We define  $P_{s,t}(x,dy)$ as the (unique) probability measure which
 makes the polynomials
$\{\QQ{n}{y}:n\in\sN\}$ orthogonal; this is possible whenever  $1+\eta\theta\geq
 0$ and  $1+x\eta\geq 0$, a condition that is
satisfied if $X_s$ has the distribution
$\pi_s(dy)=P_{0,t}(0,dy)$, see
\rf{EQ: 1+x eta}.
Since the coefficients of the three step recurrence \rf{EQ: 3-step Q}
  are bounded, it is well known that measures
$P_{s,t}(x,dy)$ have bounded support.

The next step is to prove that $P_{s,t}(x,dy)$  form a consistent
 family of measures, so that they indeed define the
transition probabilities of the Markov chain which starts at the
 origin. To this end, we need the following algebraic
relations between the
polynomials. These relations are a more complicated version of
 \cite[Theorem 1]{Bryc-Matysiak-Szablowski}
 and \cite[Lemma 3.1]{Bryc-Wesolowski-03}.

\begin{proposition}\label{THM: poly algebra}
  For  $n\geq 0$
 \begin{equation}\label{EQ:Qn}
  Q_n(z;x,u,s)=Q_{n}(y;x,t,s)+\sum_{k=0}^{n-1}B_k(y;x,t,s)Q_{n-k}(z;y,u,t)\;,
 \end{equation}
  where $B_0=1$ and
  $$
  B_1(y;x,t,s)=Q_1(y;x,t,s)-(t-s)\eta B_0\;,
  $$
  \begin{equation}\label{EQ:Bn}
  B_k(y;x,t,s)=Q_k(y;x,t,s)-t\eta B_{k-1}(y;x,t,s)\;,\;\;\;k=2,3,\dots
  \end{equation}
Additionally,  for $n\geq 1$
\begin{equation}\label{EQ:Q-p}
Q_n(y;x,t,s)=\sum_{k=0}^n
 \widetilde{B}_{n-k}(x;s)\left(p_k(y;t)-p_k(x;s)\right),
\end{equation}
where  $\widetilde{B}_k(x;s)= B_{k}(0;x,0,s)$ are linear (affine)
 functions in variable $x$.
\end{proposition}

\begin{proof}
Let
$$\phi(\zeta;y,x,t,s)=\sum_{n=0}^{\infty}\zeta^n\QQ{n}{y}$$
be the generating function of $Q_n$. Since
 $\phi(\zeta;y,x,t,s)=1+z\sum_{n=0}^{\infty}\zeta^n\QQ{n+1}{y}$,
a calculation based on recurrence \rf{EQ: 3-step Q} shows that
  $$
\phi(\zeta;y,x,t,s)=\frac{1+\zeta(t\eta+\theta-x)+\zeta^2(s+sy\eta-tx\eta+t\eta\theta )}
   {1+\zeta(t\eta+\theta-y)+\zeta^2t(1+\eta\theta)}.
     $$
From \rf{EQ:Bn} we get a similar expression for the generating
 function of $B_n$. Namely,
$$
\psi(\zeta; y,x,t,s)=\sum_{n=0}^{\infty}\zeta^nB_n(y|x,t,s)=\frac{\phi(\zeta;y,x,t,s)+\eta s \zeta}{1+\eta t\zeta}.$$
This gives
$$
\psi(\zeta; y,x,t,s)=\frac{1 + \zeta   (s \eta+\theta  -x)z+
       s  ( 1 + \eta  \theta ) \zeta^2  }
       {1 + \zeta(  t\eta  +\theta -y ) +
       t  ( 1  + \eta  \theta)  \zeta^2  }.
$$
It is now easy to verify that the two generating functions are connected by
\begin{equation}\label{EQ:87}
\phi(\zeta; z, x, u, s) - \phi(\zeta; y, x, t, s) =
 \psi(\zeta; y, x, t, s)(\phi(\zeta; z, y, u, t) - 1),
\end{equation}
which implies \rf{EQ:Qn}.
Since $\psi(\zeta,y,x,t,s)\psi(\zeta,x,y,s,t)=1$ from \rf{EQ:87} we get
$$
\phi(\zeta; z, y, u, t)=1+\psi(\zeta;x,y,s,t)(\phi(\zeta; z, x, u, s) - \phi(\zeta; y, x, t, s)).
$$
Since $p_n(x,t)=Q_n(x;0,t,0)$ setting   $x=0,s=0$ proves \rf{EQ:Q-p}.
\end{proof}
We now follow the argument from \cite[Proposition
 3.2]{Bryc-Wesolowski-03} and verify that  probability measures
$P_{s,t}(x,dy)$
are the transition probabilities of a Markov process.
\begin{proposition}\label{THM: MMM}
If $0\leq s<t<u$ and $1+\eta\theta\geq 0$, then
$$P_{s,u}(x,\cdot)=\int P_{t,u}(y,\cdot)P_{s,t}(x,dy).$$
\end{proposition}
\begin{proof}
Let
$\nu(A)=\int P_{t,u}(y,A)P_{s,t}(x,dy)$. To show that
$\nu(dz)=P_{s,u}(x,dz)$, we verify that the polynomials $Q_n(z;x,u,s)$ are
orthogonal with respect to $\nu(dz)$.
Polynomials $Q_n$ satisfy
the three-step recurrence \rf{EQ: 3-step Q}; it suffices therefore
 to show that
for $n\geq1$ these polynomials integrate to zero. Since
$
\int Q_n(z;y,u,t)
P_{t,u}(y,dz)=0$ for $k\geq 1$, by \rf{EQ:Qn} we
have
$$\int Q_n(z;x,u,s) \nu(dz)=\int  Q_{n}(y|x,t,s)P_{s,t}(x,dy)$$
$$+
\sum_{k=0}^{n-1}  \int B_k(y;x,t,s)\left(\int
Q_{n-k}(z;y,u,t)P_{t,u}(y,dz)\right)P_{s,t}(x,dy)=0.
$$

\end{proof}
For $1+\eta\theta\geq 0$, let  $(X_t)$ be the Markov process with the
 transition
probabilities  $P_{s,t}(x,dy)$, $X_0=0$.

\begin{lemma}\label{LEMMA: LM2} For $t>s, n\in\sN$ we have
\begin{equation}\label{EQ: proj}
E(p_{n}(X_t;t)|\calF_{s})=p_{n}(X_s;s).
\end{equation}
\end{lemma}
\begin{proof}  By definition, for $n\geq 1$ we have
$E(Q_n(X_t;X_s,t,s)|X_s)=0$.
Since $p_1(x,t)=x$, and $Q_1(y;x,t,s)=y-x$, by the Markov property \rf{EQ: proj} holds true for $n=1$.

Suppose that \rf{EQ: proj} holds true
for all $n\leq N$. Then \rf{EQ:Q-p}
implies
$$0=E(Q_{N+1}(X_t;X_s,t,s)|X_s)=\widetilde{B}_0(X_s;s)\left(E(p_{N+1}(X_t;t)|X_s)-p_{N+1}(X_s;s)\right).$$
Since $\tilde{B}_0=1$, this proves that $E(p_{N+1}(X_t;t)|X_s)=p_{N+1}(X_s;s)$, which by the Markov property
implies \rf{EQ: proj} for $N+1$.

\end{proof}

\begin{theorem}\label{THM: LM4} Suppose  $1+\eta\theta\geq 0$ and $(X_t)$ is the Markov process with transition
probabilities $P_{s,t}(x,dy)$,
and $X_0=0$.
Then \rf{EQ: cov}, \rf{EQ: LR},
\and  \rf{EQ: 0-Var} hold true.
\end{theorem}
\begin{proof}
Condition \rf{EQ: cov} holds true as $E(X_t)=\int
p_1(x;t)p_0(x;t)\pi_t(dx)=0$, and for $s<t$ from \rf{EQ: proj} we
get $E(X_sX_t)=E(X_sE(p_1(X_t;t)|\calF_s))=\int
p_1^2(x;s)\pi_s(dx)=\int(p_2(x;s)+(s\eta+\theta)
p_1(x;s)+s)\pi_s(dx)=s$.

Since $X_t$ are bounded,  polynomials are dense in $L_2(X_s,X_u)$.
Thus by the Markov property to prove \rf{EQ: LR} we
only need to verify
that
\begin{eqnarray}\label{EQ: LR-LR}
&E\left(p_n(X_s;s)X_tp_m(X_u;u)\right)&\\ \nonumber
=&\aaa(t|s,u) E\left(X_sp_n(X_s;s)p_m(X_u;u)\right)+
\bbb(t|s,u)E\left(p_n(X_s;s)X_up_m(X_u;u)\right)&
\end{eqnarray}
for all $m,n\in\sN$ and $0<s<t$.


For the proof of \rf{EQ: 0-Var}, we need to verify  that for any $n,m\geq 1$ and $0<s<t$
\begin{eqnarray} \label{EQ: QV-QV}
&E\left(p_n(X_s,s)X_t^2p_m(X_u,u)\right)& \\
\nonumber
= &\AAA E\left(X_s^2p_n(X_s,s)p_m(X_u,u)\right) +\BBB
E\left(X_sp_n(X_s,s)X_up_m(X_u,u)\right)& \\ \nonumber
&+\CCC E\left(p_n(X_s,s)X_u^2p_m(X_u,u)\right)
+ \alpha E\left(X_sp_n(X_s,s)p_m(X_u,u)\right)&\\ \nonumber
&+
\beta E\left(p_n(X_s,s)X_up_m(X_u,u)\right) +\DDD E\left(p_n(X_s,s)p_m(X_u,u)\right),
&\end{eqnarray}
where $\AAA,\BBB,\CCC,\DDD,\alpha,\beta$ are given by equations
\rf{EQ: A}, \rf{EQ: B}, \rf{EQ: C}, \rf{EQ: D}, \rf{EQ: alpha}, \rf{EQ: beta}:
$$
\AAA=\frac{u - t}{u - s}, \BBB=
\frac{\left( t - s \right)
         \,\left( u - t \right) }
       {\left( u - s \right) \,u}
     , \CCC=
     \frac{\left( t - s
       \right) \,t}{\left( u - s
       \right) \,u},\DDD=
  \frac{\left( t - s \right) \,
     \left( u - t \right) }{u},$$
       $$\alpha=
  \frac{\left( t - s \right)
         \,\left( u - t \right) \,
       \left( u\,\eta  - \theta
         \right) }{\left( u - s
         \right) \,u}  ,\beta=
  \frac{\left( t - s \right) \,
     \left( t - u \right) \,
     \left( s\,\eta  - \theta  \right)
       }{\left( u - s \right) \,u}
$$

It is convenient to introduce the notation
  $Ep_m^2$ for $E(p_m^2(X_s;s))$.
Recall that \rf{EQ: pn rec} implies $Ep_1^2=s$, and for $n\geq 1$
\begin{equation}\label{EQ: @}
Ep_{n+1}^2=s(1+\eta\theta)Ep_{n}^2,
\end{equation}
see  \cite[page 19]{Chihara}.

An efficient way to verify \rf{EQ: LR-LR} and \rf{EQ: QV-QV} is to use generating functions.
For $s\leq u$, let $$
\phi_0(z_1,z_2,s)=\sum_{m,n=0}^\infty  z_1^nz_2^m
E\left(p_n(X_s;s)p_m(X_u;u)\right).$$
From \rf{EQ: proj} it follows that $\phi_0(z_1,z_2,s)$ does not depend on $u$, and from \rf{EQ: @} it follows that
$$
\phi_0(z_1,z_2,s)=\frac{1-z_1z_2 \eta\theta s}{1-z_1z_2 s(1+\eta\theta)}.
$$
Consider now the generating function
$$\phi_1(z_1,z_2,s,t)=\sum_{m,n=0}^\infty  z_1^nz_2^m E\left(p_n(X_s;s)X_tp_m(X_u;u)\right).$$
From \rf{EQ: proj}  and  \rf{EQ: pn rec} we get
$$
\phi_1(z_1,z_2,s,t)=\sum_{n=0}^\infty z_1^nE\left(p_n(X_s;s)\left(X_t+z_2X_tp_1(X_t;t)+\sum_{m=2}^\infty z_2^m X_t p_m(X_t;t)\right)\right)
$$
$$
=\sum_{n=0}^\infty z_1^nE\left(p_n\left(p_1+z_2(p_2+(t\eta+\theta)p_1+tp_0)+
\sum_{m=2}^\infty z_2^m (p_{m+1}+(t\eta+\theta)p_m+t(1+\eta\theta)p_{m-1})\right)\right).
$$
Thus
$$
\phi_1(z_1,z_2,s,t)=
\left(\frac{1}{z_2}+t\eta+\theta\right)\left(\phi_0(z_1,z_2,s)-1\right)+z_2 t(1+\eta\theta)\phi_0(z_1,z_2,s)-\eta\theta t z_2,
$$
which gives
$$
\phi_1(z_1,z_2,s,t)=\frac{sz_1+tz_2+sz_1z_2(t\eta+\theta)}{1-sz_1z_2(1+\eta\theta)}.
$$
Since a calculation verifies that
$$
\phi_1(z_1,z_2,s,t)=\aaa(t|s,u)\phi_1(z_1,z_2,s,s)+\bbb(t|s,u) \phi_1(z_1,z_2,s,u),
$$
(see \rf{EQ: a+b}) from this \rf{EQ: LR-LR} follows.
Finally, for $s\leq t_1\leq t_2\leq u$ consider the generating function
$$
\phi_2(z_1,z_2,s,t_1,t_2)=\sum_{m,n=0}^\infty  z_1^nz_2^m
E\left(p_n(X_s;s)X_{t_1}X_{t_2}p_m(X_u;u)\right).$$
Another calculation based on \rf{EQ: proj}  and  \rf{EQ: pn rec} gives
$$
\phi_2(z_1,z_2,s,t_1,t_2)$$
$$=\left(\frac{1}{z_2}+t_2\eta+\theta\right)\left(\phi_1(z_1,z_2,s,t_1)-\phi_1(z_1,0,s,t_1)\right)
+z_2t_2(1+\eta\theta)\phi_1(z_1,z_2,s,t_1)-z_1z_2s\eta\theta.$$
A computer assisted calculation now verifies that
$$
\phi_2(z_1,z_2,s,t,t)$$
$$=\AAA \phi_2(z_1,z_2,s,s,s)+\BBB \phi_2(z_1,z_2,s,s,u) +\CCC \phi_2(z_1,z_2,s,u,u)+\DDD \phi_0(z_1,z_2,s)
$$
$$+
\alpha \phi_1(z_1,z_2,s,s)+\beta\phi_1(z_1,z_2,s,u),
$$
which proves \rf{EQ: QV-QV}.
\end{proof}

\section{Uniqueness}\label{Sect:Char}

We first state the main result of this section.
\begin{theorem}\label{THM: Uniqueness}
Suppose $(X_t)_{t\geq 0}$ is a centered square-integrable
separable stochastic process with covariance \rf{EQ: cov}. If $(X_t)$ satisfies \rf{EQ: LR} and \rf{EQ: 0-Var}
with $1+\eta\theta\geq 0$, then $X_t$ is the Markov process, as defined in Theorem \ref{THM: LM4}.
\end{theorem}

The proof of Theorem \ref{THM: Uniqueness} is based on the method of
 moments.
\begin{lemma}\label{Lemma:moments}
Under the assumptions of Theorem  \ref{THM: Uniqueness}
 $E(|X_t|^p)<\infty$ for all $p>0$.
\end{lemma}
\begin{proof} This result follows from \cite[Corollary 4]{Bryc85c}. To use this
 result, fix $t_1<t_2$ and let $\xi_1=t_1^{-1/2}X_{t_1}$,
$\xi_2=t_2^{-1/2}X_{t_2}$. Then their correlation
$\rho=E(\xi_1\xi_2)=\sqrt{t_1/t_2}\in(0,1)$. It
remains to notice that $E(\xi_i|\xi_j)=\rho \xi_j$ and the
variances $\Var(\xi_i|\xi_j)=1-\rho^2+a_j \xi_j$; these relations
 follow from taking the limits
$s\to 0$ or $u\to\infty$ in \rf{EQ: LR} and \rf{EQ: 0-Var}. Thus by \cite[Corollary 4]{Bryc85c}, $E(|\xi_1|^p)<\infty$ for all $p>0$.
\end{proof}
The next result is closely related to  \cite[Proposition
 3.1]{Bryc87c} and \cite[Theorem 2]{Wesolowski93}.

\begin{lemma}\label{Lemma Z} Suppose $X_t$ has  covariance \rf{EQ:
 cov}, and satisfies conditions
\rf{EQ: LR} and \rf{EQ: 0-Var}. If $0\leq s<t$ then
$E(X_t^k|\calF_s)$ is a monic polynomial of degree  $k$ in
variable $X_s$ with uniquely determined coefficients.
\end{lemma}
\begin{proof}   By Lemma \ref{Lemma:moments},
$E(|X_t^n|)<\infty$ for all $n$. Clearly, $E(X_t^k|\calF_s)$ is a
unique monic polynomial of degree $k$ when $k=0,1$. Suppose that
the conclusion holds true for  all $s<t$ and all $k\leq n$ for some integer $n\geq 1$. Multiplying \rf{EQ: LR} by $X_u^n$ and applying to both
sides  conditional expectation $E(\cdot|\calF_s)$, we get
$$
E(X_tE(X_u^n|\calF_t)|\calF_s)=\aaa X_s E(X_u^n|\calF_s)+\bbb
E(X_u^{n+1}|\calF_t).
$$
Using the induction assumption, we can write this equation as
\begin{equation}\label{EQ: ZZ1}
E(X_t^{n+1}|\calF_s)=\aaa X_s^{n+1}+  \bbb E(X_u^{n+1}|\calF_t)
+f_{n}(X_s),
\end{equation}
where $f_{n}$ is a unique polynomial of degree at most $n$.

Multiplying \rf{EQ: QV} by $X_u^{n-1}$ and applying $E(\cdot|\calF_s)$ to both sides, we get
$$
E(X_t^2E(X_u^{n-1}|\calF_t)|\calF_s)=\AAA X_s^2
E(X_u^{n-1}|\calF_s)+\BBB X_s E(X_u^{n}|\calF_t)+\CCC
E(X_u^{n+1}|\calF_t)+\dots.
$$
Using the induction assumption, we can write this equation as
\begin{equation}\label{EQ: ZZ2}
E(X_t^{n+1}|\calF_s)=(\AAA +\BBB) X_s^{n+1}+  \CCC
E(X_u^{n+1}|\calF_t) +g_{n}(X_s),
\end{equation}
where $g_{n}$ is a unique polynomial of degree   at most $n$. Since
 $\bbb-\CCC\ne 0$, subtracting \rf{EQ: ZZ1} from
 \rf{EQ: ZZ2} we get
$$
E(X_u^{n+1}|\calF_t)=\frac{\aaa-\AAA-\BBB}{\CCC-\bbb}X_s^{n+1}+h_n(X_s),
$$
where $h_n$ is a (unique) polynomial  of degree  at most $n$.

From \rf{EQ: A}, \rf{EQ: B}, \rf{EQ: C} we get
$$\frac{\aaa-\AAA-\BBB}{\CCC-\bbb}=\frac{1+\sigma u}{1+\sigma s}=1,$$  as
 $\sigma=0$. Thus $E(X_t^{n+1}|\calF_s)=X_s^{n+1}+h_n(X_s)$  is a
 monic polynomial of degree  $n+1$ in variable $X_s$
  with uniquely determined coefficients.
\end{proof}
\begin{proof}[Proof of Theorem \ref{THM: Uniqueness}] Denote by
 $(Y_t)$ the Markov process from Theorem \ref{THM: LM4}. Recall that
$Y_t$ are bounded random variables for any $t>0$ . We will show that
 by the method of moments that $(X_t)$ and $(Y_t)$ have the same
 finite dimensional distributions.

 By  Theorem \ref{THM: LM4}, process
$(Y_t)$ satisfies the assumptions of  Lemma \ref{Lemma Z}. Therefore,
 for $n\geq 0$
\begin{eqnarray}
E(Y_t^n|\calF_s)&=&Y_s^n+h_{n-1}(Y_s) \label{EQ: Z3},\\
E(X_t^n|\calF_s)&=&X_s^n+h_{n-1}(X_s) \label{EQ: Z4}
\end{eqnarray} with the same polynomial $h_{n-1}$.
From this,  we use induction to deduce that  all mixed moments are
equal. Taking $s=0$, from \rf{EQ: Z3}  and \rf{EQ: Z4} we see that
 $E(X_t^n)=E(Y_t^n)$ for all $n\in\sN,
t>0$. Suppose that for some $k\geq 1$ and all
$0<t_1<t_2<\dots<t_k$, all $n_1,\dots,n_k\in\sN$ we have
$$
E(X_{t_1}^{n_1}X_{t_2}^{n_2}\dots X_{t_k}^{n_k})=E(Y_{t_1}^{n_1}Y_{t_2}^{n_2}\dots Y_{t_k}^{n_k}).
$$
Then from \rf{EQ: Z3}  and  \rf{EQ: Z4}, by the induction assumption  we get  for
 any $t>t_k$ and $n\in \sN$
$$E\left(X_{t_1}^{n_1}X_{t_2}^{n_2}\dots X_{t_k}^{n_k}X_t^n\right)=
E\left(X_{t_1}^{n_1}X_{t_2}^{n_2}\dots
X_{t_k}^{n_k}E(X_t^n|\calF_{t_k})\right)$$
$$=
E\left(X_{t_1}^{n_1}X_{t_2}^{n_2}\dots X_{t_{k-1}}^{n_{k-1}}
X_{t_k}^{n_k}(X_{t_k}^{n}+h_{n-1}(X_{t_k}))\right)
$$
$$
E\left(Y_{t_1}^{n_1}Y_{t_2}^{n_2}\dots Y_{t_{k-1}}^{n_{k-1}}
Y_{t_k}^{n_k}(Y_{t_k}^{n}+h_{n-1}(Y_{t_k}))\right)
$$
$$=
E\left(Y_{t_1}^{n_1}Y_{t_2}^{n_2}\dots
Y_{t_k}^{n_k}E(Y_t^n|\calF_{t_k})\right)
=E\left(Y_{t_1}^{n_1}Y_{t_2}^{n_2}\dots
Y_{t_k}^{n_k}Y_t^n\right).$$ Since $t>t_k$ and $n\in\sN$ are arbitrary, this shows that
 all mixed moments of the
$k+1$-dimensional distributions match.
\end{proof}
\begin{corollary} Suppose $(X_t)$ is a Markov process from Theorem
 \ref{THM: LM4} with parameters $\eta=
\theta$. Then the process $(tX_{1/t})_{t>0}$ has the same finite
 dimensional distributions as process $(X_t)_{t>0}$.
\end{corollary}
\begin{proof} It is well known that $\rf{EQ: cov}$, and hence \rf{EQ:
 LR}, are preserved by the transformation $(X_t)\mapsto (tX_{1/t})$.
A calculation shows that if $\eta=\theta$ then the conditional
 variance
 \rf{EQ: 0-Var} is also preserved by this transformation. Thus by
 Theorem \ref{THM: Uniqueness}, both processes have the same distribution.
\end{proof}
\begin{remark}
With  more work and suitable additional  assumptions,
 Theorem \ref{THM: LM4} and  Theorem \ref{THM: Uniqueness} can
 perhaps be extended to conditional variances
 \rf{EQ: q-Var} with $\tau\ne 0$ as long as $\sigma=0, q=0$.
Generalizations to  $-1<q<1$, are hampered by the lack of suitable
 identities for the corresponding orthogonal polynomials.
When $\sigma\ne0$, an additional difficulty  arises from the fact
 that the martingale polynomial
 property \rf{EQ: proj} fails.
\end{remark}
\section{Generalized convolutions}\label{Sect: gen con}
Let $\tilde{\pi}_t$ be the measure determined by polynomials \rf{EQ:
 pn rec} with $\eta=0$.
Then  $\tilde{\pi}_t$  is a univariate distribution of the Markov
 process  $Y_t$ from
\cite[Theorem 3.5]{Bryc-Wesolowski-03} with  $\tau=q=0$.
Since this is a classical version of the free centered Poisson
 process,  it is known
that $\tilde{\pi}_t$ form a  semigroup with respect to the
 free-convolution,
 $\tilde{\pi}_{t+s}=\tilde{\pi}_t\boxplus\pi_s$.

It is somewhat surprising that there is a generalization of the
convolution that  works in a more general case; this generalization, the $c$-convolution,
is defined in
 \cite{Bozejko-Speicher91} and  studied in
 \cite{Bozejko-Leinert-Speicher}, \cite{Bozejko-Wysoczanski01}, \cite{Krystek-Yoshida03}, \cite{Krystek-Yoshida04}.

For our purposes the most convenient definition of the $c$-convolution
 is   analytic approach from \cite[Theorem 5.2]{Bozejko-Leinert-Speicher}.
According to this result, the $c$-convolution $(\mu_1,\nu_1)\boxtimes(\mu_2,\nu_2) $ is a
 binary operation  on the pairs
 of probability measures $(\mu_j,\nu_j)$, defined as follows.
Let $g_j, G_j$ be the Cauchy transforms
$$g_j(z)=\int\frac{1}{z-x}\mu_j(dx), \; G_j(z)=\int\frac{1}{z-x}\nu_j(dx)$$
On the first component of a pair, the generalized convolution acts
 just via the free convolution.
Let $k_j(z)$ be the inverse function of $g_j(z)$ in a neighborhood of $\infty$, and define
 $r_j(z)=k_j(z)-1/z$.
The free convolution $\mu$ of measures $\mu_1,\mu_2$ is defined as
 the unique probability measure
with the Cauchy transform $g(z)$ which solves the equation
$$g(z)=\frac{1}{z-r_1(g(z))-r_2(g(z))},$$
see
\cite{Voiculescu86}.

To define the second component of the $c$-convolution, let
$$R_j(z)=k_j(z)-1/G_j(k_j(z)).$$
The second component of the $c$-convolution is defined as the unique
 probability measure  $\nu$ with the Cauchy transform
$$
G(z)=\frac{1}{z-R_1(g(z))-R_2(g(z))}.
$$
We write $$
(\mu,\nu)=(\mu_1,\nu_1)\boxtimes(\mu_2,\nu_2);
$$
thus we require  that the
pair of functions $(r,R)$ as defined above be additive with respect
 to the $c$-convolution. Functions $r,R$  are the so called $r/R$-transforms
and define the $c$-free cumulants, which have interesting
combinational interpretation.

Denote by $\calL(X)$ the distribution of a random variable $X$. Let
 $Y_t$ be the free Poisson process, i.e. the Markov process from Theorem \ref{THM: LM4} with parameter $\eta=0,\theta\in\sR$.
 Let $X_t$ be
the Markov process from Theorem \ref{THM: LM4} with parameters $\eta,\theta\in \sR, 1+\eta\theta\geq 0$.

\begin{proposition}\label{THM: convolution} If $\theta=1$, then pairs of measures
 $(\calL(Y_t+t(1+\eta)),\calL(X_t+t))$ form a semigroup with respect to
the $c$-convolution,
$$(\calL(Y_{t+s}+(t+s)(1+\eta)),\calL(X_{t+s}+t+s))$$
$$=(\calL(Y_t+t(1+\eta)),\calL(X_t+t))\boxtimes (\calL(Y_s+s(1+\eta)),\calL(X_s+s)).$$
\end{proposition}
\begin{proof}
A calculation shows that $r_t(z)=\frac{t(1+\eta)}{1-z}$. Since
 $r_{t+s}(z)=r_t(z)+r_s(z)$, this verifies that indeed measures
$\calL(Y_t+t(1+\eta))$ form a semigroup with respect to free convolution.

Another calculation shows that
$$
R_t(z)=\frac{t}{1-z}.
$$
Since $R_{t+s}(z)=R_t(z)+R_s(z)$, this verifies the $c$-convolution
 property for the second component.
\end{proof}

 Measures $\pi_t$ for $\theta=1$ occur also in the Poisson Limit
 theorem for $c$-convolutions;
the Cauchy transform derived in \cite[page
 380]{Bozejko-Leinert-Speicher}
 up to centering is equivalent to \rf{EQ: G(z)}.
The conversion is accomplished by shifting argument in  \rf{EQ: G(z)}
   and making in the resulting expression
 $G(z-t)$ one of the following substitutions  $$
 \{{{\theta }\rightarrow
     1},{{\eta }\rightarrow
     {\frac{-\alpha  + \beta }
       {\alpha }}},
   {t\rightarrow {\alpha }}\}
   ,$$
   $$\{ {{\theta }\rightarrow      {-\alpha  + \beta }},
   {{\eta }\rightarrow
     {\frac{1}{\alpha }}},
   {t\rightarrow {\alpha }}\}.
$$
(The second substitution is equivalent to the first one applied to the time-reversal $tX_{1/t}$ of the bi-Poisson process.)
\begin{remark} After the first draft of this paper was written, we learned about another version of the generalized convolution,
the $\bf t$-convolution from
\cite{Krystek-Yoshida04}; this  convolution  acts on single probability measures rather than on pairs,
and could have been used in
Proposition \ref{THM: convolution} instead of the $c$-convolution.
(The case $\theta\ne 1$ still poses a challenge.)
\end{remark}
\subsection*{Acknowledgement}
We would like to thank M. Bo\.zejko and to A. Krystek for information about
the generalized convolutions and related Fock space constructions.

\bibliographystyle{siam} 
\bibliography{bipois,Vita,Wesol,q-reg}

\end{document}